\documentclass{amsart}
\usepackage{hyperref}

\usepackage{amsrefs,amsthm,amsmath,amsfonts,amssymb,graphicx,latexsym,color}

\newtheorem{theorem}{Theorem}

\theoremstyle{definition}

\newcommand{\R}{{\mathbb R}}

\newcommand{\nsm}{|}

\newcommand{\la}{\langle}
\newcommand{\ra}{\rangle}

\newcommand{\Ndim}{3}
\newcommand{\NNdim}{n}
\newcommand{\threed}{{\mathbb R}^{\Ndim}}
\newcommand{\eqdef}{\overset{\mbox{\tiny{def}}}{=}}

\begin{document}

\title{A non-local inequality and global existence}

\author[P. T. Gressman]{Philip T. Gressman}
\address{(PTG) University of Pennsylvania, Department of Mathematics, David Rittenhouse Lab, 209 South 33rd Street, Philadelphia, PA 19104-6395, USA} 
\email{gressman at math.upenn.edu}
\urladdr{http://www.math.upenn.edu/~gressman/}
\thanks{P.T.G. was partially supported by the NSF grant DMS-1101393, and an Alfred P. Sloan Foundation Research Fellowship.  }

\author[J. Krieger]{Joachim Krieger}
\address{(JK) B‰timent des MathŽmatiques
Station 8
CH-1015 Lausanne}
\email{joachim.krieger at epfl.ch}
\urladdr{http://pde.epfl.ch}
\thanks{J.K. was partially supported by SNF-grant 200021-137524.}

\author[R. M. Strain]{Robert M. Strain}
\address{(RMS) 
University of Pennsylvania, Department of Mathematics, David RittenhouseLab, 209 South 33rd Street, Philadelphia, PA 19104-6395, USA} 
\email{strain at math.upenn.edu}
\urladdr{http://www.math.upenn.edu/~strain/}
\thanks{R.M.S. was partially supported by the NSF grant DMS-0901463, and an Alfred P. Sloan Foundation Research Fellowship.}






\begin{abstract}
In this article we prove a collection of new non-linear and non-local integral inequalities.  As an example for $u\ge 0$ and $p\in (0,\infty)$ we obtain
$$
\int_{\threed} dx ~ u^{p+1}(x)
\le 
\left(\frac{p+1}{p}\right)^2
\int_{\threed} dx ~ \{(-\triangle)^{-1} u(x) \}  \nsm \nabla u^{\frac{p}{2}}(x)\nsm^2. 
$$
We use these inequalities to deduce global existence of solutions to a non-local heat equation with a quadratic non-linearity for large radial monotonic positive initial conditions.  Specifically, we improve  \cite{ksLM} to include all $\alpha\in (0, \frac{74}{75})$.
\end{abstract}

\maketitle

\thispagestyle{empty}




\section{Introduction}

In this article we study the following non-local quadratically non-linear heat equation for $\alpha >0$ given as follows:
\begin{equation}\label{eqn:model}
\partial_t u = \left\{(-\triangle)^{-1}u\right\} \triangle u + \alpha u^2,\quad  u(0, x) = u_0(x) \ge 0,
\end{equation}
where as usual  
$$
(-\triangle)^{-1}u 
=
\left(\frac{1}{4\pi |\cdot |} * u\right)(x)
=
\frac{1}{4\pi}\int_{\threed}dy~ \frac{u(y)}{|x-y|}.
$$  
We also define, for simplicity, the kernel of the Laplacian ``$(-\triangle)$'' as $G(x) = \frac{1}{4\pi |x |}$.  Note that then \eqref{eqn:model} satisfies the following conservation law:
$$
\int_{\R^3} dx ~ u(t,x) + (1-\alpha) \int_0^t ds\int_{\R^3} dx ~ |u(s,x)|^2 = \int_{\R^3} dx ~ u_0(x).
$$
In this model the variables are  $(t, x)\in [0,\infty)\times\R^{3}$. 

Equation \eqref{eqn:model} was, as far as we know, first introduced in \cite{ksLM} as a model problem for the spatially homogeneous Landau equation from 1936 \cite{MR684990}, which takes the form 
$$
\partial_t f = \mathcal{Q}(f,f).
$$
We define $\partial_i = \frac{\partial}{\partial v_i}$, and then we have the Landau collision operator
$$
\mathcal{Q}(g,f)
\eqdef 
\partial_i \int_{\threed} dv_* ~ a^{ij}(v-v_*) \left\{g(v_*) (\partial_j f)(v)  - f(v) (\partial_j g)(v_*)  \right\}.
$$
Here the projection matrix is given by
\begin{equation}
a^{ij}(v) = \frac{L}{8\pi } \frac{1}{|v|} \left( \delta_{ij} - \frac{v_i v_j}{|v|^2} \right), \quad L > 0,
\quad v = (v_1, v_2, v_3) \in \mathbb{R}^{\Ndim}.
\label{landauKernel}
\end{equation}
Throughout this article we use the Einstein convention of implicitly summing over repeated indices so that, for example,  
$a^{ij}(v)v_iv_j = \sum_{i,j=1}^\Ndim a^{ij}(v)v_iv_j$.  Above furthermore $\delta_{ij}$ is the standard Kronecker delta.
Then the following equivalent formulation of the Landau equation is well known
\begin{equation}\label{Landau}
\partial_t f = \left( a^{ij} * f \right)\partial_{i} \partial_{j} f+ L f^2,\quad (t, x)\in \R_{\geq 0}\times\threed.  
\end{equation}
See for example \cite[Page 170, Eq. (257)]{MR1942465}.  We can set $L=1$ for simplicity.  Standard references on the Landau equation include \cite{MR1942465,MR2100057,MR684990,MR1946444,MR1737547,MR2502525} and the references therein.

It is known that non-negative solutions to \eqref{Landau} preserve the $L^1$ mass.  This grants the point of view that \eqref{eqn:model} with $\alpha =1$ may be a good model for solutions to the Landau equation \eqref{Landau}. Note that \eqref{eqn:model} preserves the ``Coulomb'' singularity in \eqref{Landau}, although it removes the projection matrix.  Furthermore \eqref{eqn:model} maintains the quadratic non-linearity in \eqref{Landau}.  It appears that neither existence of global strong solutions for general large data, nor formation of singularities is known for either \eqref{Landau}, or \eqref{eqn:model} with $\alpha =1$. 
More comparisons can be found in \cite{ksLM}.

The main result of  \cite{ksLM} was to prove the following theorem.

\begin{theorem}\cite{ksLM}.  \label{thm:MainOLD} Let $0\leq\alpha<\frac{2}{3}$.  Suppose that $u_0(x)$ is positive, radial, and non-increasing with $u_0\in L^1(\threed)\cap L^{2+\delta}(\threed)$ for some small $\delta>0$.  Additionally suppose that  $-\triangle \tilde{u}_0\in L^2(\threed)$, where $\tilde{u}_0 \eqdef 
\la x \ra^{\frac{1}{2}}u_0$ and $\la x \ra \eqdef \sqrt{1+|x|^2}$. 
Then there exists a non-negative global solution, $u(t,x)$, with
\[
u(t, x)\in C^0([0, \infty), L^1\cap L^{2+\delta}(\threed))\cap C^0({\R}_{\geq 0}, H^2(\threed)),
\]
\[
\la x \ra^{\frac{1}{2}}(-\triangle)u(t, x)\in C^0([0, \infty), L^2(\threed)).
\]
Additionally the solution satisfying all of these conditions is unique.  Furthermore this solution decays toward zero at $t=+\infty$, in the following sense: 
\[
\lim_{t\to\infty}\| u(t, \cdot)\|_{L^q(\threed)}=0, \quad \forall q\in (1, (2 \wedge 1/\alpha)].
\]
Above we use the notation $(2 \wedge 1/\alpha) \eqdef\min\{2,1/\alpha\}$.  
\end{theorem}

The purpose of the present article is to improve the previous Theorem \ref{thm:MainOLD} to a substantially larger range of $\alpha\in (0, \frac{74}{75})$ in the following main theorem.

\begin{theorem}\label{thm:MainNEW} Let $0\leq\alpha<\frac{74}{75}$, and $u_0$ be as in Theorem \ref{thm:MainOLD}. Then there exists a global solution in the same spaces as in Theorem \ref{thm:MainOLD}; this solution  further satisfies
\[
\lim_{t\to\infty}\|u(t, \cdot)\|_{L^q(\threed)}=0,\quad q\in (1,75/74]. 
\]
\end{theorem}

Now our Theorem \ref{thm:MainNEW} is in some sense a consequence of the following non-linear and non-local inequality, which as far as we know is completely new:
\begin{equation}\label{modelINEQproove}
\int_{\threed} dx ~ u^{p+1}(x)
\le 
\left(\frac{p+1}{p}\right)^2
\int_{\threed} dx ~ \{(-\triangle)^{-1} u(x) \}  \nsm \nabla u^{\frac{p}{2}}(x)\nsm^2. 
\end{equation}
This inequality will hold for $p\in (0,\infty)$ and for suitable functions $u(x)\ge 0$.

Furthermore a similar inequality holds in the case of the Landau equation \eqref{landauKernel} and \eqref{Landau}.
In this situation we can prove under the same conditions the inequality
\begin{equation}\label{landauINEQproove}
\int_{\threed} dv ~ f^{p+1}
\le 
\left(\frac{p+1}{p}\right)^2
\int_{\threed} dv ~  \left( a^{ij} * f \right) \partial_i f^{\frac{p}{2}} \partial_j f^{\frac{p}{2}}. 
\end{equation}
To be precise, the inequality we prove in Theorem \ref{mainTHM} below is more general than both \eqref{modelINEQproove} and \eqref{landauINEQproove}.  These inequalities  may also be interesting on their own.

The rest of this article is organized as follows.
In the next Section \ref{lpCOMP} we supply some computations which help in proving the propagation of $L^p$ norms for solutions to \eqref{eqn:model}.  
After that in Section \ref{ldCONST} we show that for \eqref{landauINEQproove} the constant is sharp when $p=1$.  
Then in Section \ref{sec:ineq} we will state the main non-linear and non-local inequality in Theorem \ref{mainTHM}, and give its proof.
We finish the article with Section \ref{sec:implications} where we use the new inequality \eqref{modelINEQproove} and arguments from \cite{ksLM} to establish Theorem \ref{thm:MainNEW}.

\subsection{Computations regarding $L^p$ norms.}\label{lpCOMP}
 Considering the non-local equation \eqref{eqn:model}, for $p\in (0,\infty)$, we use the equation \eqref{eqn:model} obtain the following 
\begin{multline}\label{diffEQneeded}
\frac{1}{p}\frac{d}{dt} \int_{\threed} dx ~ u^{p}
=
\int_{\threed} dx ~ \partial_t u~ u^{p-1}
\\
=
-\int_{\threed} dx ~\left\{  ( G * u) \nabla u \cdot  \nabla u^{p-1} 
+\left( (\nabla G * u ) \cdot \nabla u \right) u^{p-1} 
-
\alpha  u^{p+1} 
  \right\}
\\
=
-
\frac{4}{p}
\left(\frac{p-1}{p}\right)
\int_{\threed} dx ~  ( G * u) \left| \nabla u^{\frac{p}{2}} \right|^2
-
\frac{1}{p}\int_{\threed} dx ~ u^{p+1}
+
\alpha \int_{\threed} dx ~ u^{p+1}.
\end{multline}
In other words, in order to propagate an $L^p$ norm of a solution, it would suffice for suitable non-negative functions $u$ to establish the following inequality
\begin{equation}\label{modelINEQ}
\left(\alpha - \frac{1}{p}\right)\int_{\threed} dx ~ u^{p+1} 
\le 
\frac{4}{p}
\left(\frac{p-1}{p}\right)
\int_{\threed} dx ~  ( G * u) \left| \nabla u^{\frac{p}{2}} \right|^2.
\end{equation}
In Section \ref{sec:implications} we will relate \eqref{modelINEQ} to \eqref{modelINEQproove} and prove Theorem \ref{thm:MainNEW}.  In Section \ref{sec:ineq} we will prove a more general collection of inequalities than \eqref{modelINEQ}.  Next we look at \eqref{landauINEQproove}.

\subsection{The constant in the Landau inequality.}\label{ldCONST}
In this sub-section we consider inequality \eqref{landauINEQproove}, and we argue that the constant is sharp when $p=1$.  It is well known that the Landau equation \eqref{Landau} has steady states given by the Maxwellian equilibrium, e.g. \cite{MR1942465}; for example $\mu(v) \eqdef (2\pi)^{-3/2} e^{-|v|^2 / 2}$. Then from \eqref{Landau}:
$$
0=\left( a^{ij} * \mu \right)\partial_{i} \partial_{j} \mu+  \mu^2.
$$
Similar to \eqref{diffEQneeded}, we multiply this by $\mu^{p-1}$ and integrate over $v\in \R^\Ndim$ to obtain 
$$
0=
-
\frac{4}{p}
\left(\frac{p-1}{p}\right)\int_{\R^\Ndim} dv ~ \left( a^{ij} * \mu \right)\partial_{i}\mu^{\frac{p}{2}} \partial_{j} \mu^{\frac{p}{2}}
+
\left(\frac{p-1}{p}\right)\int_{\R^\Ndim} dv ~
  \mu^{p+1}.
$$
For say $p>1$ we multiply both sides by $\frac{p}{p-1}$ and take the limit as $p\downarrow 1$ to observe
$$
4
\int_{\R^\Ndim} dv ~ \left( a^{ij} * \mu \right)\partial_{i}\mu^{\frac{1}{2}} \partial_{j} \mu^{\frac{1}{2}}
=
\int_{\R^\Ndim} dv ~
  \mu^{2}.
$$
Hence the constant in \eqref{landauINEQproove} is sharp when $p=1$.

\section{The non-linear and non-local inequality}\label{sec:ineq}

In this section we suppose that ${\bf b}(v) = (b^{ij}(v))$ is a $\NNdim \times \NNdim$ matrix for every $v\in \R^{\NNdim}$ where $i, j \in \{1, \ldots, \NNdim\}$ and say $\NNdim \ge 2$.   We furthermore suppose
\begin{equation}\label{prop1}
{\bf b}(v) = {\bf b}(-v), \quad
b^{ij}(v) \xi_i \xi_j \ge 0, \quad 
\forall  v=(v_1, \ldots, v_{\NNdim}), ~ \xi =(\xi_1, \ldots, \xi_{\NNdim}) \in  \R^{\NNdim}.
\end{equation}
Additionally we suppose that the sum over all the second derivatives in $v$ of the matrix $-{\bf b}(v-v_*)$ is a standard delta function at the point $v-v_* = 0$; precisely
\begin{equation}\label{prop2}
 -\partial_i \partial_j b^{ij}(v-v_*) = \delta_0(v-v_*).
\end{equation}
Under these basic assumptions, we have the following main inequality.

\begin{theorem}\label{mainTHM}
For suitable $g\ge 0$, under \eqref{prop1} and \eqref{prop2} we have the inequality:
\begin{equation}\notag
\int_{\R^{\NNdim}} dv ~ g^{p+1}
\le 
\left(\frac{p+1}{p}\right)^2
\int_{\R^{\NNdim}} dv ~  (b^{ij} * g) \partial_i g^{\frac{p}{2}} \partial_j g^{\frac{p}{2}} . 
\end{equation}
In this inequality we can allow any $p \in (0,\infty)$.
\end{theorem}

For simplicity and without loss of generality in the above theorem we can suppose that $g$ is a function in the Schwartz class.

Note that (for example in dimension $\NNdim = \Ndim$) inequality \eqref{modelINEQproove} follows from Theorem \ref{mainTHM} since $b^{ij}(x) = \frac{\delta_{ij}}{4\pi |x|}$ satisfies both \eqref{prop1} and \eqref{prop2}.  Similarly inequality \eqref{landauINEQproove} follows since $b^{ij}(v) =a^{ij}(v)$ from \eqref{landauKernel} is known to also satisfy \eqref{prop1} and \eqref{prop2}.  

\begin{proof}  In the proof below we can assume without loss of generality that $g$ is strictly positive; then the estimate for $g \ge 0$ will follow by approximation.
As a result of \eqref{prop1}, we consider the following quadratic form
$$
\nsm \partial g (v)\nsm_{{\bf b}(v-v_*)}^2 \eqdef  
b^{ij}(v-v_*)  (\partial_i g) (v) (\partial_j g) (v).
$$
Then for $p\in (0,\infty)$ we expand the upper bound in Theorem \ref{mainTHM} as
\begin{multline}\notag
\int_{\R^{\NNdim}} dv ~  (b^{ij} * g)(v) \partial_i g^{\frac{p}{2}} \partial_j g^{\frac{p}{2}}
=
\int_{\R^{\NNdim}} dv 
\int_{\R^{\NNdim}} dv_* ~
 b^{ij}(v-v_*)  g(v_*) (\partial_i g^{\frac{p}{2}})(v) (\partial_j g^{\frac{p}{2}})(v)
 \\
  =
\int_{\R^{\NNdim}} dv 
\int_{\R^{\NNdim}} dv_* ~
g(v_*) 
 \nsm \partial g^{\frac{p}{2}}(v) \nsm_{{\bf b}(v-v_*)}^2.
\end{multline}
By symmetry the above is 
\begin{multline}\notag
 =
 \frac{1}{2}
\int_{\R^{\NNdim}} dv 
\int_{\R^{\NNdim}} dv_* ~
g(v_*) 
 \nsm \partial g^{\frac{p}{2}}(v) \nsm_{{\bf b}(v-v_*)}^2
 +
  \frac{1}{2}
\int_{\R^{\NNdim}} dv 
\int_{\R^{\NNdim}} dv_* ~
g(v)  \nsm \partial g^{\frac{p}{2}}(v_*) \nsm_{{\bf b}(v-v_*)}^2
  \\
\ge 
\int_{\R^{\NNdim}} dv 
\int_{\R^{\NNdim}} dv_* ~
\sqrt{  g(v_*) g(v)} \nsm \partial g^{\frac{p}{2}}(v_*) \nsm_{{\bf b}(v-v_*)} \nsm \partial g^{\frac{p}{2}}(v) \nsm_{{\bf b}(v-v_*)}
   \\
=
\left( \frac{p}{p+1} \right)^2
\int_{\R^{\NNdim}} dv 
\int_{\R^{\NNdim}} dv_* ~
 \nsm \partial g^{\frac{p+1}{2}}(v_*) \nsm_{{\bf b}(v-v_*)} \nsm \partial g^{\frac{p+1}{2}}(v) \nsm_{{\bf b}(v-v_*)}.
\end{multline}
Now we use the Cauchy-Schwartz inequality for quadratic forms to obtain 
\begin{multline} \notag
\left( \frac{p}{p+1} \right)^2
\int_{\R^{\NNdim}} dv 
\int_{\R^{\NNdim}} dv_* ~
 \nsm \partial g^{\frac{p+1}{2}}(v_*) \nsm_{{\bf b}(v-v_*)} \nsm \partial g^{\frac{p+1}{2}}(v) \nsm_{{\bf b}(v-v_*)}
 \\
 \ge
 \left( \frac{p}{p+1} \right)^2
\int_{\R^{\NNdim}} dv 
\int_{\R^{\NNdim}} dv_* ~
 b^{ij}(v-v_*)   (\partial_i g^{\frac{p+1}{2}})(v) (\partial_j g^{\frac{p+1}{2}})(v_*)
    \\
=
\left( \frac{p}{p+1} \right)^2
\int_{\R^{\NNdim}} dv  ~
g^{p+1}(v). 
\end{multline}
This last computation used both \eqref{prop1}, \eqref{prop2} and standard integration by parts.
Collecting the above estimates then yields Theorem \ref{mainTHM}.
\end{proof}

In the next section we will discuss how this inequality can be used in combination with the techniques from \cite{ksLM} to obtain Theorem \ref{thm:MainNEW}.

\section{The Implications}\label{sec:implications}

In this last section, we will explain the proof of Theorem \ref{thm:MainNEW}. This proof will be derived from the developments in \cite{ksLM} and the inequality \eqref{modelINEQproove}.  In particular Theorem  \ref{thm:MainNEW} (except for the decay to zero of the $L^q$ norm) follows directly from the arguments in the proof of \cite[Theorem 1.3]{ksLM} once we show the monotonicity estimate
\begin{equation}\label{ineqPdo}
\|u(t)\|_{L^{\frac{3}{2}+\gamma}}\le \|u_0\|_{L^{\frac{3}{2}+\gamma}},
\quad 0\le t < T,
\end{equation}
for some $\gamma>0$ sufficiently small.  Here $u(t, x) \ge 0$ is a local in time solution of \eqref{eqn:model} which  is defined on $[0, T)\times\threed$ and satisfies the properties in Theorem \ref{thm:MainOLD} on $[0, T)$.

To establish \eqref{ineqPdo}, from \eqref{diffEQneeded}, \eqref{modelINEQ} and \eqref{modelINEQproove}  it suffices to check that $\alpha >0$ satisfies
$$
-
\frac{4}{p}
\left(\frac{p-1}{p}\right)+\left(\alpha - \frac{1}{p}\right)\left(\frac{p+1}{p}\right)^2 \le 0,
$$
for some $p= \frac{3}{2} + \gamma$.  This inequality is equivalent to the following 
$$
\alpha  \le 4\frac{\left(p-1\right)}{\left(p+1\right)^2}+\frac{1}{p} \eqdef h(p).
$$
Now $h(p)$ is continuous and $h(3/2) = \frac{74}{75}$ so that for any $\alpha \in (0, \frac{74}{75})$ we can find a small $\gamma >0$ such that \eqref{ineqPdo} holds.  The decay of the $L^q$ norms in given in Theorem \ref{thm:MainNEW} then follows directly from the arguments in the decay at infinity proof of \cite[Theorem 1.1]{ksLM} combined with the monotonicity in \eqref{ineqPdo}.  
\hfill {\bf Q. E. D.}

\end{document}